\documentclass{article}

\usepackage{arxiv}

\usepackage[utf8]{inputenc} 
\usepackage[T1]{fontenc}    
\usepackage{hyperref}       
\usepackage{url}            
\usepackage{booktabs}       
\usepackage{amsfonts}       
\usepackage{nicefrac}       
\usepackage{microtype}      
\usepackage{lipsum}
\usepackage{graphicx}
\graphicspath{ {./images/} }

\title{Necessary Condition for Self-organization in the BML Model with Stochastic Direction Choice
}

\author{
Marina V. Yashina  \\
  Department of Higher Mathematics\\
  Moscow Automobile and Road Construction\\
  State Technical University (MADI) \\
  Moscow, Leningradsky avenue, 64, Russia  \\
  \texttt{mv.yashina@madi.ru} \\ 
\And
 Alexander G. Tatashev  \\
  Department of Higher Mathematics\\
  Moscow Automobile and Road Construction\\
  State Technical University (MADI) \\
  Moscow, Leningradsky avenue, 64, Russia  \\
  \texttt{a-tatashev@yandex.ru} \\
}

\begin{document}
\maketitle
\begin{abstract}
 
A dynamical system is considered such that, in this system, particles move on a toroidal lattice of the dimension $N_1\times N_2$ according to a version of the rule of particle movement  in  Biham--Middleton--Levine traffic model. Particles of the first type move along rows, and the particles of the second type move along columns. A  particle can change its type. The probability that, at a step, a particle changes the type equals $q.$ It has been  proved that, if whether $q>0$ or $q=0$ and there are at least one particle of the first type and there is at least one particle of the second type, then a  necessary condition for a state of free movement to exist is that the greatest divisor of $N_1$ and $N_2$ be not less than~3. The system is represented as an algebraic structure over a set of matrices corresponding to states of the system. The theorems are formulated in terms of the algebraic structures and terms of the dynamical systems. The model is also presented as a Buslaev net. The spectrum of the net $2\times 2$ has been found.

\end{abstract}


\section{Introduction}
\hskip 18pt  In [1], a mathematical traffic model (BML model) has been introduced. In this model, particles move in a two-dimensional toroidal lattice with $N_1$ rows and $N_2$ columns. Particles move at discrete moments. There are two types of particles. Particles of the first type (red particles) move along a row from the left to the right. Particles of the second type (blue particles) move upward along a column. At any moment, a cell is vacant or occupied by a particle. According to the movement rule, at any step, a particle moves forward onto one cell if the cell ahead is vacant. The particle does not move if the cell ahead is occupied. At any step, particles of the first type move first, and the particles of the second type move later.  

In [1]--[7] some results of simulation and analytical study regarding versions of BML model were presented. In [2], an $n$-dimensional version of BML model has been introduced. Some results of simulation were analyzed. In [3], a sufficient condition for a BML model BM in a two-dimensional lattice $N\times N$ with periodic boundary conditions to result in a state of free movement (from a moment, any particle moves without delays at present and in the future) for any initial state of the system (self-organization) has been obtained. This sufficient condition is that $m<N/2,$ where $m$ is the number of particles. In [3], a sufficient condition for jam (from a moment, no particle moves) has also been obtained. In [5], it have been proved that, if $m=N/2,$ then the system also results in a state of free movement. In [7], the sufficient condition for self-organization obtained in [3] has been generilized for a version of BML model with stochastic choice of particle movement direction in $n$-dimensional lattice $N_1\times N_2\dots \times N_n$ with periodic boundary conditions. This condition is that the number of particles be not greater than the half the greatest common divisor of the numbers $N_1,N_2\dots N_n.$ 

In [7], it has been proved that a BML model may be presented as a Buslaev net. The class of dynamical systems called Buslaev nets has been obtained by A.P.~Buslaev [8], [9].    
   
In [10], a system containing two closed one-dimensional lattices called circuit. The lengths of circuits are different. The system belongs to a class of dynamical systems called Buslaev nets.  In [10], it is supposed that there is a common cell of the circuits, and particles move along circuitr according to a rule analogous to the rule
of movement in BML~model. After passing a node, with some probability  
a particle passes to the other circuit. It has been proved that, a sufficient condition for self-organization is that the number of particles be not greater than the greatest common divisor of numbers $N_1,$ $N_2,$ where $N_1$ is the number of cells in the first circuit, $N_2$ is the number of cells in the second circuit.  Under the additional condition that the probability of transition from the first circuit to the second circuit and the probability of transition from the second circuit to the first circuit are not equal to~0, this condition is also a necessary condition for self-organization.   

In this paper, we consider a system containing a toroidal lattice $N_1\times N_2.$ A particle of the first type moves along a row. A particle of the second type moves along a column.  Between moments $t$ and $t+1,$ particle changes its type with probability $0\le q< 1.$ It has been that a necessary condition for the system to result in a state of free movement for any initial state is that the greatest common divisor of the numbers $N_1,$ $N_2$ be not less than~3 under the assumption that $q>0,$ or there are at least one particle of the first type, and  at least one particle of the second type. In the proof of this statement, the approaches are used that were used in [8]. Some known facts of the number theory [11] were also used. In the deterministic version, i.~e., in the case of $q=0,$ the system is represented as a universal algebra [12] with a unary operation on matrices characterizing states of the system. This operation corresponds to the change of the system state.  In stochastic version, the system is represented by a structure with a binary relation [12] on the set of matrices of the considered type. An ordered pair of matrices $<X_1,X_2>$ belongs to this relation if the system passes in a step from the state corresponding to the matrix $X_1$ to the state corresponding to the state corresponding to matrix $X_2.$ Theorems formulated in terms of the dynamical system are also formulated and proved in terms of algebraic structures.

We give also description of the system in terms of Buslaev nets. The rows and the columns of the lattice are referred as circuits.  Each cell corresponds to a node that is a common point of two circuits. This allows to use approaches developed to study circuit networks for analysis of the BML model. An important concept is the spectrum of Buslaev net. The spectrum is the set of limits cycles (closed trajectories in the state space of the net) and the corresponding velocities of particles moving in circuits.    

We also consider as an example, in details, the $2\times 2$ system --- the simplest version of BML~model. The system is represented as a Markov
chain [13]. The spectrum of the net has been obtained.

\section{
Description of model 
}
\label{section:BuslN}
\hskip 18pt The system contains a toroidal lattice of dimension $N_1\times N_2.$ There are $m<N_1N_2$ particles. At each moment, any particle belongs to one of two types. In general case, the type of a particle may be changed.  Particles move at discrete moments, $t=0,1,2,\dots$ At each moment of time, particles of the first type move first, and then particles of the second type move. After that, type of each particle changes with probability $q$, $0\le q< 1.$ If, at time $t,$ a particle of the first type is in the cell $(i,j)$ and the cell $(i,j+1)$ (addition modulo $N_2)$ is vacant, then, at  time $t+1,$ the particle will be in the cell $(i,j+1);$ if a particle of the second type is in the cell $(i,j),$ and the cell $(i+1,j)$ (addition modulo $N_1)$ is vacant, then, at time $t+1,$ the particle will be in the cell $(i+1,j),$ $i=0,1,\dots,N_1,$ $j=0,1,\dots,N_2.$ Suppose, at time $t,$ a particle of the first type is in the cell $(i,j),$ and the cell $(i,j+1)$ (addition modulo $N_2)$ is occupied, then, at time $t+1,$ the particle will stay in the cell $(i,j).$ Suppose, at time $t,$ a particle of the second type is in the cell $(i,j),$ and the cell $(i+1,j)$ (addition modulo $N_1)$ is occupied, then, at time $t+1,$ the particle will stay in the cell $(i,j).$

\section{
Necessary condition for state of free movement
}
\label{section:BuslN}

\hskip18pt Suppose one of two following conditions holds: 1) $q>0;$ 2) $q=0,$ and there are at least one particle of the first type and at least one particle of the second type.

We use facts of the theory of linear equations in integers [11]. 
Consider the equation
$$ax + by + c = 0,\eqno(1)$$
where $a,$ $b$ are integers not equal to 0 and $c$ is an integer. There exist integer solutions of Equation (1) if and only if $c$ is divisible by the greatest common divisor of $a$ and $b.$  Let $x = x_0,$ $y = y_0$ be a solution of Equation (1). In this case, all integer solutions of Equation (1) are computed by the formulas $x = x_0 - bt,$ $y = y_0 + at,$ $t = 0, \pm 1, \pm 2,\dots$
\vskip 5pt
{\bf Theorem 1.} {\it If the greatest common divisor of $N_1$ and $N_2$ is less than 3,
and there is both at least one particle of the first type and at least one particle of the second type,  then there is no state of free movement.}
\vskip 5pt
{\bf Proof.} Suppose, at time $t_0,$ the cell $(i_1,j_1)$ is occupied by a particle of the first type, and the cell  $(i_2,j_2)$ is occupied by a particle of the second type.
Suppose $c_1=j_2-j_1$ (subtraction modulo $N_2),$ $c_2=i_2-i_1$                  (subtraction modulo $N_1),$ $c=c_2-c_1.$

If non-negative integers $x$ and $y$ satisfy the equation
$$N_2x-N_1y=c,\eqno(2)$$
then, at time $t_0+N_2x,$ a delay of the second type particle occurs if, during the time interval $(t_0,t_0+N_2x)$ delays do not occur and the types of particles do not change.   

If non-negative integers $x$ and $y$ satisfy the equation 
$$N_2x-N_1y=c+1,\eqno(3)$$
then, at time $t_0+N_2x,$ a delay of the first type particle occurs if, during the time interval  $(t_0,t_0+N_2x)$ delays do not occur and the types of particles do not change. 

If the greatest common divisor $N_1$ and $N_2$ is not greater than 2, then the right hand of at least one of Equations (2), (3) is divisible by the greatest common divisor, and therefore, for at least one of Equations (2) and (3), there are integers  $x$ and $y$ satisfying the equation.   

Taking into account that, with non-negative  probability, particles do not change the sign at a finite time interval, we get that, at time $t_0,$ the system is not in a state of free movement. This completes the proof of the theorem.
\vskip 5pt
The condition that the greatest common divisor of $N_1$ and $N_2$ is equal to 3 is not a sufficient condition for that a state of free movement does not exist.     
\vskip 5pt
{\bf Example 1.}  Suppose $N_1=N_2=3,$ $q=0.$ There are only one particle of the first type and one type of the second type. At time $t=0,$ the particle of the first type is in the cell $(0,0),$ and the particle of the second type is in the cell  $(0,2).$ Then the particles move without delays, and, at time $t=3,$  the particles will be in the same cells as at time $t=0.$ Thus there exists the state of free movement, and the greatest common divisor of  $N_1$ and $N_2$ equals 3.

\section{
Algebra corresponding to deterministic version of dynamical system
}
\label{section:BuslN}
\hskip 18pt Consider an algebra $<A,f>,$ where $A=A(N_1,N_2,m_1,m_2)$ is a set of matrices $N_1\times N_2$ containing $m_1$ elements equal to 1, $m_2$ elements equal to 2 and the other elements to 0; $f$ is a unary operation that takes $X=(x_{ij})\in A$ to $Y=(y_{ij})\in A.$ Let us define the operation~$f.$ 

First the matrix $U(X)=(u_{ij})\in A$ is obtained. If $x_{ij}=0$ and $x_{i,j-1}=0$ (subtraction modulo $N_2),$ then $u_{ij}=0,$ $i=1,\dots,N_1,$ $j=1,\dots,N_2.$ If $x_{ij}=0$ and $x_{i,j-1}=1$ (subtraction modulo $N_2),$ then $u_{ij}=1.$ If $x_{ij}=0$ and $x_{i,j-1}=2$ (subtraction modulo $N_2),$ then $u_{ij}=0.$ If $x_{ij}=1$ and $x_{i,j+1}=0$ 
(addition modulo $N_2),$ then $u_{ij}=0.$ If $x_{ij}=1$ and $x_{i,j+1}=1$ (addition modulo $N_2),$ then  $u_{ij}=1.$ If $x_{ij}=1$ and $x_{i,j+1}=2,$ then $u_{ij}=1.$ 
If $x_{ij}=2,$ then $u_{ij}=2.$ 

Let us define the matrix $Y(U).$ The $U(X)$ is taken to $Y$ as following. If $u_{ij}=0$ and $u_{i-1,j}=0,$ then $y_{ij}=0,$ $i=1,\dots,N_1,$ $j=1,\dots,N_2.$ If $u_{ij}=0$ and $u_{i-1,j}=1,$ then $y_{ij}=0.$ If $u_{ij}=0$ and $u_{i+1,j}=2,$ then $y_{ij}=2.$ If $u_{ij}=1,$ then $y_{ij}=1.$ If $u_{ij}=2$ and $u_{i+1,j}=0,$ then $y_{ij}=0.$ If $u_{ij}=2$ and $u_{i+1,j}=1,$ then $y_{ij}=2.$ If $u_{ij}=2$ and $u_{i+1,j}=2,$ then $u_{ij}=2.$ 

Suppose 
$$f(X)=Y(X)=Y(U(X)).\eqno(4)$$

If the matrix $X$ corresponds to the state of the system at time $t$ before the movement of particles of the first type,  
the matrix $U(X)$ corresponds to the state of the system at time $t$ after the movement of particles of the first type, 
the matrix $Y(X)$ corresponds to the state of the system at time $t+1$ before the movement of particles of the first type,

We say that $X=(x_{ij})\in A$ belongs to $B$ if, for any $i,$ $j$ such that $x_{ij}=1,$ the condition $x_{i,j+1}=0$ (addition modulo~$N_2)$ holds. 
We say that $U=(u_{ij})\in B$  if, for any $i,$ $j$ such that, if $u_{ij}=1,$ the condition $u_{i+1,j}=0$ (addition modulo~$N_1)$ holds,  $i=1,\dots,N_1,$ $j=1,\dots,N_2.$ 

The set $B$ is the set of matrices corresponding to the states of systems such that any particle moves.

Suppose $f^1(X)=f(X),$ $f^t(X)=f(f^{t-1}(X)),$ $X\in A,$ i.e. $f^t$ the $t$-fold composition of the operation $f,$ $n=2,3,\dots$  

Suppose $X(0)=(x_{ij}(0))=X,$ $X(t)=(x_{ij}(t))=f^t(X)$ $U(X(t))=(u_{ij}(t))$ is the matrix obtained from the matrix
$X(t)$ according to the above algorithm of transition from the matrix $X$ to the matrix $U(X).$

We say that the matrix $X\in A$ belongs to the set $C$ if there exists a number $t_0$ such that, for any $t\ge t_0,$   
the matrix $f^t(X)$ belongs to the set $B.$

The set $C$ is the set of matrices corresponding to the states of free movement.

Let $d$ be the greatest common divisor of  $N_1$ and $N_2.$  

\vskip 5pt
 {\bf Theorem 2.} {\it If the greastest common divisor of $N_1$ and $N_2$ is less than~3, and there are at least one particle of the first 
type and at least one particle, then no matrix belongs to set $C.$   
\vskip 5pt
Proof.} Suppose, for the matrix $X=(x_{ij}),$ $x_{i_1j_1}=1,$ $x(i_2j_2)=1.$
Suppose $c_1=j_2-j_1$ (addition modulo $N_2),$ $c_2=i_2-i_1$ (subtraction modulo $N_1),$ $c=c_2-c_1.$

Suppose non-negative integers $x$ and $y$ satisfy Equation (2) or Equation (3).

Suppose $t=t_0+N_2x.$ Then $f^t(X)\notin C$ if, for $t\le t_0+N_2x,$ $f^t(X)\notin C.$     

If the greatest common divisor of $N_1$ and $N_2$ is nor greater than 2, then the right hand (2), (3) is divisible by this number, and, therefore, for, at least one equation (4), (5), there are non-negative integers $x,$ $y$ satisfying this equation.   

Taking into account then, with non-zero probability, the type particles changes durung a finite time interval, we get that the state of the system at time $t_0$ is not a state of free movement. Theorem~2 has been proved.

\section{
Relational structure with stochastic version of dynamical system
}
\hskip 18pt Consider the structure $<D,r>,$  where $D=D(N_1,N_2,m)$ is a set of matrices of dimension $N_1\times N_2$ such that $m$ elements of matrices  are equal to~1, and the other elements are equal to~0;  $r$ is a binary relation. We shall give the definition of this relation. 

An ordered pair $<X,Y>$ belongs to the relation $r$ (we write $XrY)$ if the matrix $Y$ can be obtained from the matrix $X,$ if, first, the values of $m_2$ arbitrarily chosen elements of the matrix $X,$ equal to~1, change to 2 and pass to the matrix $Y$ according to the rule described in Section~6. 

Let $r^t$ be a binary relation defined as follows. An ordered pair of matrices $<X,Y^t>$ belongs to the relation $r^t$ if there  exists a matrix $Y^{t-1}$ such that $Xr^{t-1}Y^{t-1}$ and $Y^{t-1}rY,$  $t=2,3,\dots,$ $Y^0=X.$
\vskip 5pt
If $X(t)$ is the matrix corresponding to the state of the system before the moving of particles of the first type, and  then $XrU,$ then the matrix $U$ is a matrix 
such that, the system, with a positive probability, will be in the state corresponding to matrix~$U$ after moving particles of the first type.  

If $X$ is the matrix corresponding to the state of the system before the movement of particles of the first type, and $XrU,$ then the matrix $U$ is a matrix 
such that, the system, with a positive probability, will be in the state corresponding to matrix~$U$ after moving particles of the first type.  

If $U$ is the matrix corresponding to the state of the system after the movement of particles of the first type, and $UrY,$ then the matrix $Y$ is a matrix 
such that, the system, with a positive probability, will be in the state corresponding to matrix~$Y$ at time $t+1$ before the movement of particles of the first type.  
\vskip 5pt
The definitions of sets $B$ and $C$ are analogous to definitions given in Section~4. The set $B$ is the set of matrices corresponding to states such that any particle moves with probability~1. The set $C$ is the set of matrices corresponding to states of free movement.
  
\vskip 5pt
 {\bf Theorem 3.} {\it If the greatest common divisor of the numbers $N_1$ and $N_2$ is less than~3,  then the set $C$ does not contain any matrix. 
}
\vskip 5pt
Theorems 2, 3 are analogous to Theorem 1.

From Theorems 2, 3, it follows that the statements of 
Theorems 1 may be generalized for a more general system than the 
system studied in Sections~2--4. At any step, the probability that the type of a particle changes may depend on the index of the particle,
the type of the particle, etc.    

\section{
$2\times 2$ system
}
  
\hskip 18pt Suppose $m=n=2,$ Fig. 1.

\begin{figure}[ht!]
\centerline{\includegraphics[width=300pt]{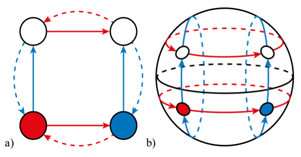}}
\caption{a) $2\times 2$ BML model as Buslaev net, b) subset of $S^2$}
\end{figure}

Consider a deterministic Markov chain containing $2\times 3^4 = 162$ states corresponding to the BML model of the dimension $2\times 2.$ This model may be interpreted as a Buslaev net.

Suppose, the matrix $X$ corresponds to a state configuration,
$$
X=\left(
\begin{array}{c c}
x_{21}&x_{22}\\
x_{11}&x_{12}\\
\end{array}
\right).
$$
Suppose 
$$(m,m_1,m_2;x_{11},x_{12},x_{21},x_{22};s)$$ 
is the state such 
that $m$ is the number of particles; $m_1$ is the number of particles of the 
first type; $m_2$ is the number of particles of the second type, $x_{11},$ $x_{12},$ $x_{21},$ $x_{22}$ are the elements of the matrix $X;$ $s$ is the type of moving particles.

Suppose, at time $t_0,$ the Markov chain is in the state $i_0.$ If there exists $t_1>t_0$ such that, at time $t_1,$ the Markov chain is in the same state $i_0.$ Then the state $i_0$ is called recurrent [13]. 

Any recurrent state of a deterministic finite Markov chain belongs to a closed trajectory in the state space of the chain. These trajectories are called cycles.     

The number of states belonging to a cycle is called the period of the cycle. 

Suppose $T$ is the period of a cycle, and $A$ is the number of transitions of a particle during the cycle. Then the value of $A/T$ is called the velocity of the particle, and the velocity of each particle corresponds to the cycle.  

The following results are obtained for the $2\times 2$ system.
\vskip 3pt
Each of 162 states is recurrent, and hence each state belongs to a cycle. 
\vskip 3pt 
There are two states $(0,0,0;0,0,0,0;1),$ $(0,0,0;0,0,0,0;2)$ with no particle.
These states belongs to a cycle with period 2.
\vskip 3pt
There are 16 states with one particle. There are 4 cycles with period 4, and the velocity~1. These cycles are
$$(1,1,0;0,0,0,1;1)\to (1,1,0;0,0,1,0;2)\to (1,1,0;0,0,1,0;1)\to (1,1,0;0,0,0,1;2)\to \dots,$$     
$$(1,1,0;0,1,0,0;1)\to (1,0,0;1,0,0,0;2)\to (1,1,0;1,0,0,0;1)\to (1,1,0;0,1,0,0;2)\to \dots,$$     
$$(1,0,1;0,0,0,2;1)\to  (1,0,1;0,0,0,2;2)\to (1,0,1;0,2,0,0;1)\to (1,0,1;0,0,0,2;2)\to \dots,$$
$$(1,0,1;0,0,2,0;1)\to  (1,0,1;0,0,2,0;2)\to (1,0,1;2,0,0,0;1)\to (1,0,1;2,0,0,0;2)\to \dots$$
\vskip 3pt
There are 8 states with two the first type particles located in different rows. These states belong to two cycles with period equal to~2 and the velocity equal to~1. These cycles are
$$(2,2,0;0,1,0,1;1)\to (2,2,0;1,0,1,0;2)\to (2,2,0;1,0,1,0,1;1)\to (2,2,0;0,1,0,1,2)\to \dots,$$
$$(2,2,0;0,1,1,0;1)\to (2,2,0;1,0,0,1;2)\to (2,2,0;1,0,0,1;1)\to (2,2,0;1,0,0,1;2)\to \dots$$
\vskip 3pt
There are 8 states with two the second type particles located in different columns. These states belong to two cycles with period equal to~4 and the velocity equal to~1. These cycles are
$$(2,0,2;0,0,0,2,2;1)\to (2,0,2;0,0,2,2;2)\to (2,0,2;2,2,0,0;1)\to (2,0,2;2,2,0,0;2)\to \dots,$$
$$(2,2,0;0,2,2,0;1)\to (2,2,0;0,2,2,0;2)\to (2,2,0;2,0,0,2;1)\to (2,2,0;2,0,0,2;2)\to \dots$$
\vskip 3pt
There are 4 states with two the first type particles located in the same row. These states belong to two cycles with period equal to~2 and the velocity 0. These cycles are
$$(2,2,0;0,0,1,1;1)\to (2,2,0;0,0,1,1;2)\to \dots,$$
$$(2,2,0;1,1,0,0;1)\to (2,2,0,1,1,0,0;2)\to \dots$$
\vskip 3pt
There are 4 states with two the second type particles located in the same row. These states belong to two cycles with period equal to~2 and the velocity 0. These cycles are
$$(2,0,2;0,2,0,2;1)\to (2,0,2;0,2,0,2;2)\to \dots,$$
$$(2,0,2;2,0,2,0;1)\to (2,0,0;2,0,2,0;2)\to \dots.$$
\vskip 3pt
There are 24 states with two particles of different types. These states belong to 4 cycles with period equal to~6 and velocity equal to~2/3. These cycles are
$$(2,1,1;0,0,1,2;1)\to (2,1,1;0,0,1,2;2)\to
(2,1,1;0,2,1,0;1)\to (2,1,1;0,2,0,1;2)\to $$
$$(2,1,1;0,2,0,1;1)\to (2,1,1;0,2,1,0;2)\to \dots $$
\vskip 3pt
$$(2,1,1;0,0,2,1;1)\to (2,1,1;0,0,2,1;2)\to
(2,1,1;2,0,0,1;1)\to (2,1,1;2,0,1,0;2)\to $$
$$(2,1,1;2,0,1,0;1)\to (2,1,1;2,0,0,1;2)\to \dots $$
\vskip 3pt
There are 8 states with three particles of the first type. These states belong to 2 cycles with period equal to~4, and the velocity of two particles is equal to~0, and the velocity of one particle is equal to~1. The cycles are 
$$(3,3,0;0,1,1,1;1)\to (3,3,0;1,0,1,1;2)\to 
(3,3,0;1,0,1,1;1)\to (3,3,0;0,1,1,1;2)\to \dots
$$
\vskip 3pt
$$(3,3,0;1,1,0,1;1)\to (3,3,0;1,1,1,0;2)\to 
(3,3,0;1,1,1,0;1)\to (3,3,0;1,1,0,1;2)\to \dots
$$
\vskip 3pt
There are 8 states with three particles of the second type. These states belong to 2 cycles with period equal to~4, and the velocity of two particles is equal to~0, and the velocity of one particle is equal to~1. These cycles are
$$(3,0,3;0,2,2,2;1)\to (3,0,3;0,2,2,2;2)\to 
(3,0,3;2,2,0,2;1)\to (3,0,3;2,2,0,2;2)\to \dots
$$
\vskip 3pt
$$(3,0,3;2,0,2,2;1)\to (3,0,3;2,0,2,2;2)\to 
(3,0,3;2,2,0,2;1)\to (3,0,3;2,2,0,2;2)\to \dots
$$
\vskip 3pt
There are 16 states with two particles of the first type in different rows and one particle of the second type. These states belong to two cycles with period equal to~8 and the velocity equal to~1/2. The cycles are 
$$(3,2,1;0,1,1,2;1)\to (3,2,1;1,0,1,2;2)\to 
(3,2,1;1,2,1,0;1)\to (3,2,1;2,2,0,1;2)\to $$
$$(3,2,1;2,2,0,1;1)\to (3,2,1;2,2,1,0;2)\to 
(3,2,1;0,1,1,2;1)\to (3,2,1;1,0,1,2;2)\to \dots
$$
\vskip 3pt
$$(3,2,1;0,1,2,1;1)\to (3,2,1;1,0,2,1;2)\to 
(3,2,1;1,0,2,1;1)\to (3,2,1;0,1,2,1;2)\to $$
$$(3,2,1;0,1,2,1;1)\to (3,2,1;1,0,2,1;2)\to
(3,2,1;1,0,2,1;1)\to (3,2,1;0,1,2,1;2)\to \dots
$$
\vskip 3pt
There are 16 states with two particles of the second type in different columns and one particle of the first type. These states belong to two cycles with period equal to~8 and the velocity equal to~1/2. The cycles are 
\vskip 3pt
$$(3,1,2;0,2,2,1;1)\to (3,2,1;0,2,2,1;2)\to 
(3,2,1;2,2,0,1;1)\to (3,2,1;2,1,2,1;2)\to $$
$$(3,2,1;0,1,2,1;1)\to (3,2,1;1,0,2,1;2)\to
(3,2,1;1,0,2,1;1)\to (3,2,1;0,1,2,1;2)\to \dots
$$
\vskip 3pt
$$(3,2,1;2,0,1,2;1)\to (3,2,1;2,0,1,2;2)\to 
(3,2,1;2,2,1,0;1)\to (3,2,1;2,2,0,1;2)\to $$
$$(3,2,1;0,2,2,1;1)\to (3,2,1;0,2,2,1;2)\to
(3,2,1;2,2,0,1;1)\to (3,2,1;2,2,1,0;2)\to \dots
$$

\vskip 3pt
There are 8 states with two particles of the first type in the same row and one particle of the second type. These states belong to 4 cycles with period equal to~2 and the velocity equal to~0. The cycles are 
$$(3,2,1;0,2,1,1;1)\to (3,2,1;0,2,1,1;2) \to \dots, $$  
$$(3,2,1;1,1,0,2;1)\to (3,2,1;1,1,0,2;2)\to \dots, $$
$$(3,2,1;1,1,2,0;1)\to (3,2,1;1,1,2,0;2)\to \dots, $$
$$(3,2,1;2,0,1,1;1)\to (3,2,1;2,0,1,1;2)\to \dots $$
\vskip 3pt
There are 8 states with two particles of the second type in the same column and one particle of the first type. These states belong to 4 cycles with period equal to~2 and the velocity equal to~0. The cycles are 
$$(3,2,1;0,2,1,2;1)\to (3,2,1;0,2,1,2;2)\to \dots, $$
$$(3,2,1;1,2,0,2;1)\to (3,2,1;0,2,2,1;2)\to \dots, $$
$$(3,2,1;2,0,2,1;1)\to (3,2,1;2,0,2,1;2)\to  \dots, $$
$$(3,2,1;2,1,2,0;1)\to (3,2,1;2,1,0,2;2)\to \dots $$
\vskip 3pt
There are 32 states with four particles. These states belong to 16 cycles with period equal to~2 and the velocity equal to~0. The cycles are 
$$(4,4,0;1,1,1,1;1)\to (4,4,0;1,1,1,1;2)\to \dots,$$
$$(4,3,1;1,1,1,2;1)\to (4,4,0;1,1,1,2;2)\to \dots,$$
$$(4,3,1;1,1,2,1;1)\to (4,4,0;1,1,2,1;2)\to \dots,$$
$$(4,3,1;1,2,1,1;1)\to (4,3,1;1,2,1,1;2)\to \dots,$$
$$(4,3,1;2,1,1,1;1)\to (4,3,1;2,1,1,1;2)\to \dots,$$
$$(4,2,2;1,1,2,2;1)\to (4,2,2;1,1,2,2;2)\to \dots,$$
$$(4,2,2;1,2,1,2;1)\to (4,2,2;1,2,1,2;2)\to \dots,$$
$$(4,2,2;1,2,2,1;1)\to (4,2,2;1,2,2,1;2)\to \dots,$$
$$(4,2,2;2,1,1,2;1)\to (4,2,2;2,1,1,2;2)\to \dots,$$
$$(4,2,2;2,1,2,1;1)\to (4,2,2;2,1,2,1;2)\to \dots,$$
$$(4,2,2;2,2,1,1;1)\to (4,2,2;2,2,1,1;2)\to \dots,$$
$$(4,1,3;1,2,2,2;1)\to (4,2,2;1,2,2,2;2)\to \dots,$$
$$(4,1,3;2,1,2,2;1)\to (4,2,2;2,1,2,2;2)\to \dots,$$
$$(4,1,3;2,2,1,2;1)\to (4,1,3;2,2,1,2;2)\to \dots,$$
$$(4,1,3;2,2,2,1;1)\to (4,4,0;1,1,2,1;2)\to \dots,$$
$$(4,0,4;2,2,2,2;1)\to (4,0,4;2,2,2,2;2)\to \dots$$
\vskip 3pt
Thus we have the following

{\it There are 16 states with one particle and the velocity of any particle is equal to~1.

There are 48 states with 2 particles. Among these states, there are 16~states with velocity of particles equal to~1, 24~states with the average velocity of particles equal to~2/3, 8~states with velocity of particles equal to~0. Therefore, under the assumption that the number of particles is equal to~2, the average velocity is equal to~2/3.    

There 64 states with three particles. Among these states, there are 
32~states with the average velocity equal to~1/2, 16~states with velocity equal to~1. Hence, under the assumption that the number of particles is equal to~2, the average velocity is equal to~2/3.

There are 16 states with one particle and the velocity of any particle is equal to~0.
}
\vskip 3pt

Thus we have found the spectrum, i.~e., the set of cycles with corresponding 
velocities.

\section{Conclusion}

\hskip 18pt Versions of BML-models are helpful to obtain optimal solutions for problems of optimal traffic control.  

 A system is studied such that, in this system, particles move
along a toroidal lattice $N_1\times N_2$ according to rules generalizing the rules of the BML model. A necessary condition for free movement has been found. According to this condition, if the greatest common divisor of  $N_1$ and $N_2$ is less than~3, and either $q>0$ or $q=0$ and there is at least one particle of the first type,  and there is at least one particle of the second type, then there exists no state of free movement. The condition is also formulated in terms of algebraic structures.

An example has been considered. We have obtained the spectrum of the Buslaev net corresponding to the considered BML-model. 

\bibliographystyle{unsrt}  


\end{document}